\documentclass[12pt]{article}

\newtheorem{theorem}{Theorem}[section]

\newtheorem{proposition}[theorem]{Proposition}
\newtheorem{corollary}[theorem]{Corollary}

\newfam\msbfam
\font\tenmsb=msbm10  scaled \magstep1 \textfont\msbfam=\tenmsb
\font\sevenmsb=msbm7 scaled \magstep1 \scriptfont\msbfam=\sevenmsb
\font\fivemsb=msbm5  scaled \magstep1 \scriptscriptfont\msbfam=\fivemsb
\def\Bbb{\fam\msbfam \tenmsb}

\def\RR{{\Bbb R}}
\def\CC{{\Bbb C}}
\def\QQ{{\Bbb Q}}
\def\NN{{\Bbb N}}
\def\ZZ{{\Bbb Z}}
\def\II{{\Bbb I}}
\def\TT{{\Bbb T}}
\def\BB{{\Bbb B}}

\def\ss{\subseteq}
\def\ra{\rightarrow}
\def\btu{\bigtriangleup}

 \def\HollowBox #1#2{{\dimen0=#1 \advance\dimen0 by -#2       
       \dimen1=#1 \advance\dimen1 by #2                       
        \vrule height #1 depth #2 width #2                    
        \vrule height 0pt depth #2 width #1                   
        \llap{\vrule height #1 depth -\dimen0 width \dimen1}%
       \hskip -#2                                             
       \vrule height #1 depth #2 width #2}}                   
 \def\BoxOpTwo{\mathord{\HollowBox{6pt}{.4pt}}\;}             

\def\endpf{\hfill $\BoxOpTwo$ \smallskip \\ }

\newfam\msbbfam
\font\tenmsbb=msbm10  scaled \magstep1 \textfont\msbbfam=\tenmsbb
\font\sevenmsbb=msbm7  scaled \magstep1 \scriptfont\msbbfam=\sevenmsbb
\font\fivemsbb=msbm5    scaled \magstep1 \scriptscriptfont\msbbfam=\fivemsbb

\usepackage{graphicx}

\usepackage{amsmath}

\begin{document}

\begin{center}
\Large \bf Homage \`{a} A. M. Gleason et I. J. Schark\footnote{{\bf Subject 
Classification Numbers:}  30H80, 30H05, 30H50, 32A35, 32A37, 32A38, 32A65.}\footnote{{\bf Key Words:}  corona problem, unit ball,
fiber, S\u{\i}lov boundary, maximal ideal space, analytic disc, uniform algebras, Banach algebras.}
\end{center}
\vspace*{.12in}

\begin{center}
Steven G. Krantz
\end{center}

\date{\today}

\begin{quote}
{\bf Abstract:}   We study the maximal ideal space of $H^\infty(B)$, where
$B$ is the unit ball of $\CC^n$.  Following the lead of Gleason and Schark,
we analyze Gleason parts, fibers, the S\u{\i}lov boundary, and other aspects of
this Banach algebra.  Our work here makes good use of the inner functions construction
of Aleksandrov and Hakim/L\o w/Sibony, particularly as formulated by Rudin.
\end{quote}
\vspace*{.25in}

\markboth{STEVEN G. KRANTZ}{GLEASON AND SCHARK}

\section{Introduction}

The corona problem on the disc was formulated by S. Kakutani in 1941.  It was solved by Lennart Carleson
in 1962 (see [CAR]).  There was not a great deal of work on the matter in the intervening 21 years.
But a few papers stand out.  Among these are [SCH] by the eponymous group I. J. Schark and
[GLE] by A. M. Gleason.

In fact it may certainly be said that [SCH] laid the foundations for all future studies
of the corona.  Here the idea of fiber was introduced, and many of its basic properties
established.  Some of the exotic structure of the maximal ideal space was proved.

Today there is great interest in studying the corona problem on the unit ball $B$ in $\CC^n$.
But little is known.  It is certainly worthwhile to re-examine some of the ideas introduced
in [SCH] and [GLE] in the new context of $B$ and to see what is true and what has changed.
That is what we endeavor to do in the present paper.

We discover, naturally, that some of the basic ideas go through in a rather pro forma fashion.
Others require new effort, or new ideas, or new formulations.  And still others fail or at least need to
be modified.   Studying the corona in
this more general context is not only of interest in its own right, but it sheds light
on the classical situation.   It also illustrates some nice applications of the inner
functions construction of [ALE], [HAS], [LOW].	See also [RUD1].

The paper [GLE] is the source of the idea of Gleason part.  Again, this is a concept that
has not seen much development in the several variable context.  We hope to initiate
such a study here.

\section{Fiber Basics}

Fix a complex Euclidean space $\CC^n$ once and for all.
Let $H^\infty(B)$ denote the algebra of bounded analytic functions
on the unit ball $B \ss \CC^n$.  Equipped with the uniform norm,
$H^\infty(B)$ is a commutative Banach algebra.  Denote its maximal ideal
space by ${\cal M}$.  As usual, we topologize ${\cal M}$ with the
weak-star topology.  Thus ${\cal M}$ is a compact, Hausdorff space, and $\widehat{H^\infty(B)}$ 
(where the accent \ \ $\widehat{}$ \ \ denotes the Gelfand transform) is
a uniformly closed subalgebra of continuous functions on ${\cal M}$.

Let $\sigma$ be the rotationally invariant area measure on $S = \partial B$.
Of course we know (see, for instance, [KRA1]) that each $f \in
H^\infty(B)$ has associated to it a $\sigma$-almost-everywhere-defined
radial boundary limit function $\widetilde{f} \in
L^\infty(\partial B)$. We sometimes use the symbol $S$ to
denote $\partial B$.  We will write $H^\infty(S)$ instead of $H^\infty(B)$ with
no danger of ambiguity.

We shall discuss below the Gelfand transform $f \mapsto
\widehat{f}$ which maps $H^\infty(B)$ to a uniformly closed
algebra of continuous functions on ${\cal M}$. Also note that
there is a natural projections $\pi$ of ${\cal M}$ onto the
closed unit ball in complex space $\CC^n$, obtained by sending
$\varphi \in {\cal M}$ to $(\varphi(z_1), \varphi(z_2), \dots,
\varphi(z_n))$. Here we think of each $z_j$ as a coordinate
function, hence an element of $H^\infty(B)$.

The map $\pi$ is injective over the open ball $B$.  To see this assertion notice that,
if $\pi(\varphi) = \pi(\varphi')$, then $\varphi(z_1, z_2, \dots, z_n) = \varphi'(z_1, z_2, \dots, z_n)$ 
so that $\varphi$ and $\varphi'$ are the same function.   So we see that the natural
injection $i$ of $B$ into ${\cal M}$, which sends $\lambda \in B$ into the point evaluation
functional at $\lambda$, is a homeomorphism of $B$ onto an open subset
$\bigtriangleup$ of ${\cal M}$.  The remaining set ${\cal M} \setminus \bigtriangleup$ is 
a closed set of multiplicative homomorphisms and is mapped by $\pi$ onto the sphere
$S = \partial B$.

The set ${\cal M} \setminus \bigtriangleup$ is decomposed by $\pi$ into disjoint closed {\it fibers}
(i.e., the inverse images of various points in $S$ under $\pi$).  More precisely, if
$\alpha \in S$, then the fiber over $\alpha$ is
$$
{\cal M}_\alpha = \biggl \{ \varphi \in {\cal M}: \pi(\varphi) = \bigl ( \varphi(z_1), \varphi(z_2), \dots, \varphi(z_n) \bigr ) = \alpha \biggr \} = \pi^{-1}(\alpha) \, .
$$
The action of the rotation group $SU(n)$ on $S$, and hence on $H^\infty(S)$, shows that
any given fiber ${\cal M}_\alpha$ is homeomorphic with any other given
fiber ${\cal M}_{\alpha'}$.  We shall see later that these fibers are
{\it not necessarily} analytically equivalent.

Later on we shall identify the S\u{\i}lov boundary for the
Banach algebra $H^\infty(B)$. The description will be as
follows. Since, as noted above, each element $f \in
H^\infty(B)$ has a boundary function $\widetilde{f} \in
L^\infty(S)$, we then have a natural, continuous map $\tau$ of
the (extremely disconnected) space of maximal ideals of
$L^\infty(S)$ into the space ${\cal M}$. It turns out that
$\tau$ is a homeomorphism, with image the S\u{\i}lov boundary
of $H^\infty(B)$.  On the one hand, the S\u{\i}lov boundary
is a subset of $\overline{\bigtriangleup} \setminus \bigtriangleup$.
On the other hand, it does {\it not} exhaust ${\cal M} \setminus \bigtriangleup$.

Later sections contain additional properties of the fibers.

\section{Algebraic Properties}

Now we formulate and prove our first concrete result about the ideas introduced
thus far.  In what follows, we let $\varphi_\lambda$ denote the 
point evaluation functional at $\lambda \in \overline{B}$.

\begin{proposition} \sl 
The projection $\pi$ defined above is a
continuous mapping of ${\cal M}$ onto the closed unit ball $B$
in space. If $\bigtriangleup = \pi^{-1}(B)$, then $\pi$ maps
the open set $\bigtriangleup$ homeomorphically onto the open
ball $B$.
\end{proposition}
{\bf Proof:}  Certainly $\pi$ is a continuous, $\CC^n$-valued function
on ${\cal M}$.	It maps ${\cal M}$ into the closed unit ball because
$\|(z_1, z_2, \dots, z_n)\| = 1$ when $(z_1, z_2, \dots, z_n) \in B$.
Each point of the open ball $B$ is in the range of $\pi$, because $\pi(\varphi_\lambda) = \lambda$.
Since $\pi({\cal M})$ is compact and contains $B$, it contains
the entire closed ball $\overline{B} = S \cup B$.

We note again, in a slightly different fashion, that $\pi$ is one-to-one over $B$.  That is to say, if $\varphi \in {\cal M}$
and $\pi(\varphi) = (\varphi(z_1), \varphi(z_2), \dots, \varphi(z_n)) = \lambda$ with
$\|\lambda\| < 1$, then $\varphi(f) = 0$ for every $f$ of the form
$f(z) = (z_1 - \lambda_1)g_1(z) + (z_2 - \lambda_2) g_2(z) + \cdots + (z_n - \lambda_n) g_n(z)$.
So $\varphi(f) = 0$ whenever $f(\lambda) = 0$.  Thus $\varphi$ must be point
evaluation $\varphi_\lambda$ at $\lambda$.

If $\bigtriangleup = \pi^{-1}(B) = \{\varphi_\lambda: \lambda \in B\}$, then $\bigtriangleup$ is
an open subset of ${\cal M}$. Either on $\bigtriangleup$ or on $B$, the topology is
the weak-$*$ topology.  Thus $\pi$ is a homeomorphism of $\bigtriangleup$ onto $B$.
\endpf
\smallskip \\

Certainly $\pi$ maps ${\cal M} \setminus \bigtriangleup$ onto the sphere $S$.   Also
$\pi$ is not one-to-one over $S$, and we have designated the fiber over $\alpha \in S$
as ${\cal M}_\alpha$.  An element of ${\cal M}_\alpha$ should be thought of as roughly
like ``evaluation at $\alpha$; that is, a homomorphism $\varphi$ of $H^\infty(B)$ which
sends each $f \in H^\infty(B)$ into a sort of limiting value of $f(\lambda)$ as $\lambda$ approaches
the boundary point $\alpha$.  We shall flesh this idea out later.  For now we concentrate
on more elementary properties of ${\cal M}_\alpha$.

The algebra $H^\infty(B)$ is of course rotationally invariant.  That is to say,
if $f \in H^\infty(B)$ and $\sigma \in SU(n)$, then the function
$z \mapsto f(\sigma z)$ is still in $H^\infty(B)$.  And the new function
has the same norm as the old.  Define the rotation
$$
(R_\sigma f)(z) \equiv f(\sigma z) \, .
$$
We may say that $R_\sigma$ is an {\it automorphism} of $H^\infty(B)$.  The
adjoint mapping $R_\sigma^*$ defined by
$$
(R_\sigma^* \varphi)(f) = \varphi (R_\sigma f) 
$$
is thus a homeomorphism of ${\cal M}$ onto ${\cal M}$.  So the rotation
group acts as a discrete group of homeomorphisms of ${\cal M}$.  By this we mean
that each $R_\sigma^*$ is an isolated point in the group.  Now we may see that
the fibers ${\cal M}_\alpha$ are homeomorphic.  For let $\alpha$, $\beta$ be distinct
points of the sphere $S$.  Suppose that $\sigma \in SU(n)$ is such that
$$
\beta = \sigma \alpha \, .
$$
Then the map $R_\sigma^*$ carries ${\cal M}_\alpha$ homeomorphically onto ${\cal M}_\beta$.  For
$(\varphi(z_1), \varphi(z_2), \dots,$ \\
$\varphi(z_n)) = \alpha$ if and only if
$(R_\sigma^* \varphi)(z) = \varphi (\sigma z) = \sigma \alpha = \beta$.

\section{More on the S\u{\i}lov Boundary}

As usual, let $L^\infty(S)$ be the space of all essentially bounded, $\sigma$-measurable, complex-valued
functions on the sphere $S$.  With the usual operations and the usual norm, $L^\infty(S)$ is a commutative
Banach algebra.  With each element $f \in H^\infty(B)$ we may associate, by way of the radial
limit process, a unique $L^\infty(S)$ function $\widetilde{f}$ (see [KRA1] for the details).
Note that $\widetilde{f}$ exists almost everywhere on $S$ and
$$
\|f\| = \|\widetilde{f}\|_\infty \, .
$$
One can recover $f$ from $\widetilde{f}$ either by way of the Poisson integral
formula or one of the other classical integral formulas in the subject (again see [KRA1]).

The mapping
$$
f \longrightarrow \widetilde{f}
$$
is an isometric isomorphism of $H^\infty(B)$ onto a closed subalgebra
of $L^\infty(S)$.  We often denote, without danger of confusion,
the latter subalgebra by $H^\infty(S)$ or just $H^\infty$.  In complex
dimension 1, we may characterize the elements of $H^\infty(S)$ in terms of
orthogonality to certain trigonometric monomials (essentially the F. and M.
Riesz theorem).  Such a characterization is not possible in the several
variable setting, although see [CHG].

Observe that $L^\infty(S)$ is not only a commutative Banach algebra under
the essential supremum norm, but it is also closed under complex conjugation
(recall the Stone-Weierstrass theorem).  In the ensuing discussion, we use $X$ 
to denote the maximal ideal space of $L^\infty(S)$.  Thus $L^\infty(S)$ is isometrically
isomorphic to the algebra $C(X)$ of all continuous, complex-valued functions on
a Hausdorff space $X$ (see [LOO, p. 88]).  By Gelfand's theory, the
space $X$ is nothing other than the space of complex multiplicative homomorphisms
(maximal ideals) of $L^\infty(S)$.  For $g \in L^\infty(S)$, we let $\widehat{g}$ denote
the corresponding continuous function on $X$ (it is of course given by the Gelfand 
transform).  

If $M$ is any measurable subset of the sphere $S$, we let $\chi_M$ be the characteristic
functions of $M$.  Since the functions $\chi_M$ generate $L^\infty$ in norm, then
the functions $\widehat{\chi_M}$ generate $C(X)$.  Because $\chi_M^2 = \chi_M$,
we see that a basic open subset of $X$ has the form
$$
\{x \in X:  \widehat{\chi_M}(x) = 0 \} \, ,
$$
where $M \subseteq S$ is measurable.  A set of this type is also obviously
closed, so we see that $X$ is totally disconnected.  It can even be shown
that $X$ is extremely disconnected.

There is a natural continuous mapping (as described above) $\tau$ of $X$ into ${\cal M}$ because
we can identify $H^\infty(B)$ with the subalgebra
$H^\infty$ of $L^\infty$.  A point $x \in X$ is a complex homomorphism of $L^\infty$, and by $\tau(x)$
we mean the complex homomorphism of $H^\infty(B)$ obtained by restricting $x$
from $L^\infty(S)$ to $H^\infty(S)$ (identified in the obvious way with $H^\infty(B)$).  

The classical proof of Theorem 4.2 below (see [SCH]) uses constructions involving 
the conjugate function.  These are unavailable in the several complex variable
context.   We shall use instead an interesting result of W. Rudin [RUD1]:

\begin{proposition} \sl
Let $\psi \in L^\infty(S)$.  Suppose that there is an $f \in H^\infty(B)$, with
$f$ never 0, such that $\psi \geq |\widetilde{f}|$ \ \  $\sigma$-almost-everywhere.  Further assume
that $\psi/|\widetilde{f}|$ is almost lower semicontinuous.  Then there is
a function $F \in H^\infty(B)$ such that $|\widetilde{F}| = \psi$ a.e. $\sigma$.
\end{proposition}

It is really quite remarkable how nicely Rudin's result fits our needs in the following theorem.

\begin{theorem} \sl
The mapping $\tau$ is a homeomorphism of $X$ {\it into} ${\cal M}$.  The range
$\Gamma \equiv \tau(X)$ is the S\u{\i}lov boundary for $H^\infty(B)$.  That is, $\Gamma$
is the smallest closed subset of ${\cal M}$ on which
every function $\widehat{f}$, with $f \in H^\infty(B)$, attains its maximum modulus.
\end{theorem}
\noindent {\bf Proof:}   Fix a point $x_0 \in X$.  Let $U$ be a basic neighborhood (as described above) of the
point $x_0$.  Then $U$ has the form
$$
U = \{x \in X: \widehat{\chi_M}(x) = 0\} \, ,
$$
where $M$ is a measurable subset of the sphere $S$ such that $\widehat{\chi_M}(x_0) = 0$.
Let us apply Proposition 4.1 with 
$$
\psi(x) =  \exp \left ( 1 - \chi_M \right ) 
$$
and $f \equiv 1$.  Then clearly $\psi/f$ is almost lower semicontinuous if we simply 
assume that $M$ is the closure of an open set with piecewise smooth boundary.  Also $\psi \geq \widetilde{f} \equiv 1$.
So we obtain a bounded holomorphic function $F$ whose boundary function satisfies
$|\widetilde{F}| = \psi$ a.e.  Plainly
$$
|F| = \left \{ \begin{array}{lcr}
                   e \ & \ \hbox{on} \ & \  U \\
		   1 \ & \ \hbox{on} \ & \ X \setminus U \, .
	       \end{array}
      \right.
$$

These simple calculations show us that the functions $\widetilde{f}$ for $f \in H^\infty(B)$ separate the points
of $X$.  Also there is plainly no proper closed subset of $X$ on which all such $\widetilde{f}$ attain their
maximum modulus.  

The map $\tau$ is easily seen to be continuous.  Since $X$ is compact and we know that $\tau$ is one-to-one, 
we have that $\tau$ is a homeomorphism.  For any $f \in H^\infty(B)$, 
$$
\sup_{\cal M} |\widetilde{f}| = \|f\| = \|\widetilde{f}\|_\infty = \sup_X |\widehat{\widetilde{f}}| \, .
$$
Since $\widehat{f}(\tau(x)) = \widehat{\widetilde{f}}(x)$, we see that $\Gamma \equiv \tau(X)$ is the smallest
closed subset of ${\cal M}$ in which each $\widehat{f}$ attains its maximum modulus.   Thus it is the S\u{\i}lov
boundary.
\endpf
\smallskip \\

We now have the following information about ${\cal M}$.  First, ${\cal M}$ contains the (image of the) open unit
ball $B$.  It also contains $\Gamma$, which is the S\u{\i}lov boundary for $H^\infty(B)$.  The set
$\Gamma$ is homeomorphic to the extremely disconnected maximal ideal space of $L^\infty(S)$.  The closure
$\overline{\btu}$ of $\btu$ contains $\Gamma$.  This is because, since
$$
\sup_\btu |\widehat{f}| = \|f\| = \sup_{\cal M} |\widehat{f}| \, ,
$$
each $\widehat{f}$ will attain its maximum value on $\overline{\btu}$.  The maximum modulus
principle of course tells us that $\Gamma$ is contained in $\overline{\btu} \setminus \btu$. 

There are also points in $\overline{\btu} \setminus \btu$ which are
{\it not} in the S\u{\i}lov boundary $\Gamma$. For consider an
inner function $g$ on the unit ball $B$ (see [KRA1, Ch.\ 9]
for the details). Of course the boundary function
$\widetilde{g}$ has modulus 1 $\sigma$-almost everywhere. So
$\widetilde{g}$ does not vanish on $\Gamma$. It is easily
arranged that $g$ not vanish on $B$ (in Aleksandrov's [ALE] construction
of inner functions this is quite explicit), but that at least one
radial boundary limit for $g$ equals 0. We may assume that this
radial limit is at the point ${\bf 1} = (1, 0,0, \dots, 0)$.
Thus $\widetilde{g}$ vanishes somewhere on $S$.
Since $g$ is nonvanishing on $B$, the zero of $\widetilde{g}$ occurs
at a point of $\overline{B} \setminus B$ that is not in $\Gamma$.

\section{Values Taken on the Fibers}

Certainly we have defined the fibers ${\cal M}_\alpha$ over
various points $\alpha$ of the unit sphere $S$. To be
specific, the fiber over $\alpha \in S$ consists of those
complex multiplicative homomorphisms $\varphi$ of
$H^\infty(B)$ such that $(\varphi(z_1), \varphi(z_2), \dots,
\varphi(z_n)) = \alpha$.  Thus it is clear that, for any $f \in H^\infty(B)$ 
which is continuously extendable to the closed ball $\overline{B}$, the
function $\widehat{f}$ is constant on each fiber ${\cal M}_\alpha$.
This is so because such an $f$ is the uniform limit of polynomials
in $z$ (by the power series expansion).   Now we can make this
assertion more precise.  Namely, the continuity of $f$ at any one boundary point
$\alpha$ of $S = \partial B$ implies that $\widehat{f}$ is constant
on ${\cal M}_\alpha$.

\begin{theorem}  \sl
Let $f \in H^\infty(B)$ and let $\alpha$ be a point of the sphere $S$.
Let $\lambda_j$ be a sequence of points in the open ball $B$ such that
\begin{enumerate}
\item $\lambda_j \ra \alpha$,
\item $\zeta = \lim_{j \ra \infty} f(\lambda_j)$ exists.
\end{enumerate}
Then there is a complex multiplicative homomorphism $\varphi$ of $H^\infty(B)$ which lies
in the fiber ${\cal M}_\alpha$ and for which $\varphi(f) = \zeta$.
\end{theorem}
{\bf Proof:}  Let $I$ be the set of all functions $g \in H^\infty(B)$ such that 
$\lim_{j \ra \infty} g(\lambda_j) = 0$.  Then $I$ is a proper, closed ideal in the
algebra $H^\infty(B)$.  Therefore $I$ is contained in some maximal
ideal $M$.  Let $\varphi$ be the complex multiplicative homomorphism of $B$ of which $M$ is
the kernel.  Then $z - \alpha$ is in $I$, and so is $f - \zeta$.  Therefore
$$
\varphi(z) = \alpha \ \ \ \hbox{and} \ \ \ \varphi(f) = \zeta \, .
$$
So $\varphi$ is the required multiplicative homomorphism.
\endpf
\smallskip \\

\begin{theorem} \sl
Let $f \in H^\infty(B)$ and $\alpha \in S$.  A necessary and sufficient
condition for $\widehat{f}$ to be constant on ${\cal M}_\alpha$ is that $f$ be
continuously extendable to $B \cup \{\alpha\}$.
\end{theorem}
{\bf Proof:}  First assume that $f$ is continuously extendable to $B \cup \{\alpha\}$.
This means that there is a complex number $\zeta$ so that $f(\lambda_j) \ra \zeta$
whenever $\lambda_j \ra \alpha$.  We claim that $\varphi(f) = \zeta$ for all $\varphi \in {\cal M}_\alpha$.
We may suppose that $\zeta = 0$.  We may also assume, after composing
the picture with a unitary rotation, that $\alpha = (1, 0, 0, \dots, 0)$.  Then set
$g(z_1, z_2, \dots, z_n) = ((z_1 + 1)/2, z_2/2, z_3/2, \dots, z_n/2)$.  We see that $g(\alpha) = (1,0, 0, \dots, 0)$ and
$|g| < 1$ elsewhere on $\overline{B}$.  Since $f$ is continuous at $\alpha$, if we
set $f(\alpha) = 0$, we may see that $(1 - g^j)f$ converges uniformly on compact sets to $f$
as $j \ra \infty$.  If $(\varphi(z_1), \varphi(z_2), \dots, \varphi(z_n)) = \alpha$, then
$\varphi(g) = 1$ and $\varphi(1 - g^j) = 0$, so that $\varphi(f) = 0$.

If $\widehat{f}$ has the constant value $\zeta$ on the fiber ${\cal M}_\alpha$, then the preceding theorem
implies that $f(\lambda_j) \ra \zeta$ whenever $\lambda_j \ra \alpha$.  If one defines $f(\alpha) = \zeta$,
then $f$ is continuous on $B \cup \{\alpha\}$.
\endpf
\smallskip \\ 

\begin{theorem} \sl
Let $f \in H^\infty(B)$ and $\alpha \in S$ and suppose that there is a complex,
multiplicative homomorphism $\varphi \in {\cal M}_\alpha$ such that $\varphi(f) = 0$.
Then there exists a sequence $\{\lambda_j\}$ such that $\lambda_j \ra \alpha$ and
$f(\lambda_j) \ra 0$.
\end{theorem}
{\bf Proof:}  If there is no such sequence $\{\lambda_j\}$, then there is a neighborhood
${\cal N}$ of the point $\alpha$ such that $|f(\lambda)| \geq \delta > 0$ for all $\lambda \in B \cap {\cal N}$.
We may slice the ball $B$ with a complex line $\ell$ that passes through ${\cal N}$ and also
through the origin.  Rotating coordinates, we may as well suppose that this complex
line is $\ell = \{z: z_2 = z_3 = \cdots = z_n = 0\}$.   Since there is a natural
injective map from $H^\infty(D)$ (living in the line $\ell$) to $H^\infty(B)$, then
there is a natural surjective dual map from ${\cal M}_{H^\infty(B)}$ to ${\cal M}_{H^\infty(D)}$.
So if $f$ is annihilated by $\varphi \in {\cal M}_{H^\infty(B)}$, then
$f$ restricted to $\ell$ will be annihilated by the image
of $\varphi$ in ${\cal M}_{H^\infty(D)}$.  Therefore this restricted
function satisfies the hypotheses of Theorem 4.3 in [SCH].  As a result,
we derive a contradiction as in that source and find that there is a sequence
$\{\lambda_j\}$ as claimed.
\endpf
\smallskip \\

\begin{corollary} \sl
Let $f \in H^\infty(B)$ and $\alpha \in S$. Then the range of $\widehat{f}$ on the fiber
${\cal M}_\alpha$ consists precisely of those complex numbers $\zeta$ for which there is
a sequence $\{\lambda_j\}$ in $B$ with
\begin{tabbing}
\null \qquad \qquad \qquad \qquad \= {\bf (i)} \ \ \= $\lambda_j$ \ \ \ \ \ \= $\rightarrow$ \ \= $\alpha$ \, , 
\smallskip \\
				  \> {\bf (ii)}    \> $f(\lambda_j)$        \> $\rightarrow$   \> $\zeta$ \, .
\end{tabbing}
\end{corollary}

Now we wish to comment on the topological nature of the decomposition
$$
{\cal M} \setminus \bigtriangleup = \bigcup_{\alpha \in S} {\cal M}_\alpha \, .  \eqno (5.5)
$$
We already know that the fibers are all homeomorphic under the action
of the rotation group.  Thus one might hope that (5.5) is a product decomposition,
that is, that ${\cal M} \setminus \bigtriangleup$ is naturally homeomorphic to 
the product space $S \times {\cal M}_\alpha$ for some fixed fiber ${\cal M}_\alpha$.
As a point set, ${\cal M} \setminus \bigtriangleup$ can be naturally identified
with $S \times {\cal M}_{\bf 1}$.  Here boldface {\bf 1} denotes the point
$(1, 0, 0, \dots, 0)$ in the boundary $S$ of the ball $B$. 

We associate with each $\alpha \in S$ a unitary rotation $\sigma_\alpha$ that
maps the point $\alpha \in S$ to the point ${\bf 1} = (1, 0, \dots, 0) \in S$.
If $\varphi \in {\cal M}_\alpha$, we associate
with $\varphi$ the ordered pair $(\alpha, R_{\sigma_\alpha}^* \varphi)$, where $R_{\sigma_\alpha}^*$
is the adjoint of the rotation induced by $\alpha$ (see Section 3).  Now the mapping
$$
\varphi \longmapsto (\alpha, R_{\sigma_\alpha}^* \varphi)
$$
is a one-to-one correspondence (as sets) between ${\cal M} \setminus \bigtriangleup$ and
$S \times {\cal M}_{\bf 1}$.  However, it is {{\it not} a homeomorphism.  If it were a homeomorphism,
then the orbit
$$
\{R_\beta^* \varphi \}_{\beta \in S}
$$
of any $\varphi$ under the rotation group would be a continuous cross section
of ${\cal M} \setminus \btu$ over $S$.  This is impossible because such an orbit is not a closed
subset of ${\cal M} \setminus \btu$.  In fact, there are no continuous cross-sections of ${\cal M} \setminus \btu$
over $S$ whatsoever.  To see this last assertion, let us make a number of observations about the
topological nature of the decomposition (5.5).

Let $W_+$ be the union of the fibers ${\cal M}_\alpha$ for $\hbox{Im}\, \alpha_1 > 0$ and let 
$W_{-}$ be the union of all ${\cal M}_\alpha$ with $\hbox{Im}\, \alpha_1 < 0$.  Denote the closures
of these two sets by $\overline{W}_+$ and $\overline{W}_{-}$.  Then we have:
\begin{enumerate}
\item[{\bf (i)}]   The intersection of $\overline{W}_+$ and $\overline{W}_{-}$ is empty.
\item[{\bf (ii)}]    Each point in the fiber ${\cal M}_{\bf 1}$ is in either $\overline{W}_+$ or $\overline{W}_{-}$ 
or else it is in neither.  All three cases actually occur.
\item[{\bf (iii)}]     If $\{\varphi_j\}$ is any sequence of points in ${\cal M} \setminus \btu$ which converges,
then all but a finite number of the $\varphi_j$ lie in the same fiber ${\cal M}_\alpha$.
\item[{\bf (iv)}]    If $\Sigma$ is any function from $S$ into ${\cal M} \setminus \btu$ such that $\pi \circ \Sigma$ is
the identity (so $\Sigma$ is a section of $\pi$), then $\Sigma(S)$ is not closed.  In particular, such a section $\Sigma$
cannot be continuous.
\end{enumerate}
\vspace*{.22in}

Now let us prove these four assertions.
\smallskip \\

\noindent {\bf Proof of (i):}  Use Proposition 4.1 above to construct a bounded holomorphic
function $f$ on $B$ so that $|f|$ equals 1 on $\{\hbox{Im}\, \alpha_1 > 0\}$ and equals
$e$ on $\{\hbox{Im}\, \alpha_1 < 0\}$.  Then the preceding corollary tells us that
$$
|\widehat{f}(\varphi)| = \left \{ \begin{array}{lcr}
                               1 \ & \ \hbox{if} \ & \ \varphi \in W_+ \, , \\
			       e \ & \ \hbox{if} \ & \ \varphi \in W_{-} \, .
				  \end{array}
			 \right.
$$
Thus we see that $\overline{W}_+$ and $\overline{W}_{-}$ are disjoint, just becuse the continuous function
$\widehat{f}$ satisfies $|\widehat{f}| = 1$ on $\overline{W}_+$ and
$|\widehat{f}| = e$ on $\overline{W}_{-}$.
\smallskip \\

\noindent {\bf Proof of (ii):}   This is a tautology:  obvious by inspection.
\smallskip \\


\noindent {\bf Proof of (iii):}  Let $\{\varphi_j\}$ be a sequence of points
in ${\cal M} \setminus \btu$.  If there are more than a finite number
of points of the sphere $S$ among the images $\alpha_j \equiv \pi(\varphi_j)$,
then we need to show that $\{\varphi_j\}$ does not converge.   Passing to 
a subsequence if necessary, we may suppose that the $\alpha_j$ are distinct.
Choose for each $j$ a spherical cap $A_j$ centered at $\alpha_j$.  We may assume
that these spherical caps are pairwise disjoint.  Now let us apply Proposition 4.1
to obtain an $H^\infty$ function $f_0$ on the ball $B$ so that (for $\alpha \in S$)
$$
|f_0(\alpha)| = \left \{ \begin{array}{lcr}
                    e^{(-1)^j} \ & \ \hbox{if} \ & \ \alpha \in A_j \, ,  \\
			 0     \ & \           \ & \hbox{otherwise}.
			 \end{array}
		\right.
$$
So we see that 
$$
\varphi_j(f_0) = \left \{ \begin{array}{lcr}
                        e \ & \ \hbox{if} \ & \ j \ \hbox{is even} \, ,  \\
			1/e \ & \ \hbox{if} \ & \ j \hbox{is odd} \, .
			  \end{array}
		 \right.
$$
That shows that $\varphi_j$ cannot converge.
\smallskip \\

\noindent {\bf Proof of (iv):}  This follows immediately from {\bf (iii)}.  If $\Sigma$
is a function from $S$ into ${\cal M} \setminus \btu$ such that $\pi \circ \Sigma$ is the identity,
then let $\varphi_0 = \Sigma({\bf 1})$.  Choose a sequence
of distinct points $\alpha_j$ on the sphere $S$ such that $\alpha_j \ra {\bf 1}$, and let 
$\varphi_j = \Sigma(\alpha_j)$.  If $\sigma(S)$ were closed, then $\varphi_0$ would be the only
cluster point of the sequence $\{\varphi_j\}$.  But, by {\bf (iii)}, the sequence
$\{\varphi_j\}$ cannot converge to $\varphi_0$ because the $\alpha_j$ are
distinct.
\medskip \\

Our statement {\bf (ii)} gives some idea of the bizarre nature of the topology
on ${\cal M} \setminus \btu$.  It says that, in any given fiber ${\cal M}_\alpha$,
some of the points can be approached from points in fibers ${\cal M}_\beta$ near to but distinct from
${\cal M}_\alpha$, and some points of ${\cal M}_\alpha$ cannot be approached from points in any 
other fibers ${\cal M}_\beta$.  These are mutually exclusive possibilities which can both occur.  It can be
shown that no point of the S\u{\i}lov boundary is of the second type.  It is not known whether
${\cal M}_\alpha$ has any interior.  This would be tantamount to solving the corona problem
on the ball $B$.

Two other topological problems, each having intrinsic interest, are these:  {\bf (a)}  Since ${\cal M} \setminus \btu$
contains an extremely disconnected set $\Gamma$, we may ask whether ${\cal M} \setminus \btu$ is connected;
and {\bf (b)} we may ask whether each fiber ${\cal M}_\alpha$ is connected.

\section{Embedding a Disc in a Fiber}

Now we consider complex structure in ${\cal M}$.

A mapping $\psi$ from the unit disc $D$ into ${\cal M}$ is called {\it analytic}
if $\widehat{f} \circ \psi$ is analytic on $D$ for each $f \in H^\infty(B)$. 
If $\{\psi_j\}$ is a sequence of analytic maps of $D$ into ${\cal M}$,
then the compactness of ${\cal M}$ guarantees that there is a cluster
point $\psi$ of $\{\psi_j\}$ in the space of maps of $D$ into ${\cal M}$.  
We may see that $\psi$ is analytic because, for each $f \in H^\infty(B)$, the
sequence $\widehat{f} \circ \psi_j$ is uniformly bounded and hence 
uniformly equicontinuous on each compact subset of $D$.  [Note that this reasoning
is a form of Montel's theorem and particularly of its proof.]

Now we shall construct an analytic map $\psi$ of $D$ into the space ${\cal M}$ which
is {\bf (i)} a homeomorphism and {\bf (ii)} maps $D$ into the fiber ${\cal M}_{\bf 1}$.
The ``disc'' $\psi(D)$ in the fiber ${\cal M}_{\bf 1}$ will have the property that the
restriction of $\widehat{H^\infty(B)}$ to $\psi(D)$ consists precisely of all
bounded analytic functions on this disc.

Let $L$ be the linear fractional transformation
$$
L(\lambda_1, \lambda_2, \dots, \lambda_n) = \left ( \frac{(1 + i)\lambda_1 - i}{(1 - i) + i\lambda_1},
					    \frac{(1 + i)\lambda_2/\sqrt{2}}{(1 - i) + i\lambda_1}, \dots, 
					    \frac{(1 + i)\lambda_n/\sqrt{2}}{(1 - i) + i\lambda_1} \right ) \, .
$$
One may check directly that $L$ is a biholomorphic mapping of the ball $B$ to itself, and that $L$
in fact maps the closure $\overline{B} = B \cup S$ to itself.  Notably, the point ${\bf 1} = (1, 0, 0, \dots, 0)$ 
is the only fixed point of this mapping (note that $- {\bf 1} = (-1, 0, 0, \dots, 0)$ is {\it not} a fixed
point).  If we let $L^j$ denote the composition of $L$ with itself $j$ times, then it is easy to calculate
that
$$
L^j(\lambda) = \left ( \frac{\lambda_1 + ji(\lambda_1 - 1)}{1 + ji(\lambda_1 - 1)},
			       \frac{\lambda_2}{1 + ji(\lambda_1 - 1)}, \dots, 
			       \frac{\lambda_n}{1 + ji(\lambda_1 - 1)}\right ) \, .	 \eqno (6.1)
$$

Now let $\psi_j$ be the map of $B$ into ${\cal M}$ defined by
$$
\psi_j(\lambda) = \pi^{-1} (L^{2^j} (\lambda)) \, .
$$
In other words, we see that $\psi_j$ maps $B$ into $\btu$ and $\psi_j(\lambda)$ is the complex 
homomorphism of $H^\infty(B)$ which evaluates each $f$ in $H^\infty(B)$ at the pont $L^{2^j}(\lambda)$
in $B$.  Clearly $\psi_j$ is an analytic map of $B$ into ${\cal M}$ (in fact even into $\btu$).  Let
$\psi$ be a cluster point of the sequence of maps $\{\psi_j\}$ so that $\psi$ is an analytic map
of $B$ into ${\cal M}$.  It is easy to see that $\psi$ must map $B$ into the fiber ${\cal M}_{\bf 1}$.  For, if
we fix a point $\lambda \in B$, then 
$$
\lim_{j \ra \infty} L^j (\lambda) = {\bf 1} \, .
$$
This shows that $\pi(\psi(\lambda)) = {\bf 1}$ for each $\lambda \in B$.  

Fix a unit vector $\xi \in \CC^n$.  Our claim is that $\psi$ is a one-to-one analytic mapping of $B$ into ${\cal M}_{\bf 1}$.  We set
$$
f(\lambda) = \left [  \prod_{j=0}^\infty \xi \cdot L^{-2^j} (\lambda) \right ] (\xi \cdot \lambda) \, .  \eqno (6.2)
$$

We see from line (6.1) that, on any compact subset of the ball $B$, 
$$
|L^j (\lambda) \cdot \xi - 1| \leq K \cdot \left ( \frac{1}{|j|} \right ) \, , \qquad \hbox{for} \ \ |j| > 0 \, .  \eqno (6.3)
$$
Since $|L^j| \leq 1$, we have shown that the infinite product (6.2) converges uniformly on
compact subsets of $B$ to a function $f$ in $H^\infty(B)$ with $\|f\| \leq 1$.   We shall
utilize this $f$ to show that $\psi: B \ra {\cal M}_{\bf 1}$ is a homeomorphism.

We claim that 
$$
\widehat{f}(\psi(\lambda)) = \lambda \cdot \xi \qquad \hbox{for} \ \ \lambda \in B \, .	   \eqno (6.4)
$$
For, using (6.3) for $\lambda$ in a compact subset of $B$, we see that
\begin{eqnarray*}
|f(L^{2^k}(\lambda)) - \lambda \cdot \xi | & = & \left | \xi \cdot L^{2^k}(\lambda) \prod_{j=0}^\infty \xi \cdot L^{2^k - 2^j}(\lambda) - \lambda \cdot \xi \right | \\
                                & = & |\lambda \cdot \xi| \left | \xi \cdot L^{2^k}(\lambda) \prod_{j=0}^{k-1} \xi \cdot L^{2^k - 2^j} 
                                      \cdot \prod_{j=k+1}^\infty \xi \cdot L^{2^k - 2^j}(\lambda) - 1 \right | \\
				& \leq & |\lambda| \left [ |\xi \cdot L^{2^k}(\lambda) - 1| + \sum_{j=0}^{k-1} 
                                      |\xi \cdot L^{2^k - 2^j}(\lambda) - 1| \right.  \\
                                &      &  \left. + \sum_{j = k+1}^\infty |\xi \cdot L^{2^k - 2^j}(\lambda) - 1| \right ] \, , \\
\end{eqnarray*}
where we have used in the last line some standard estimates for infinite products.  This last is
\begin{eqnarray*}
 & \leq & |\lambda| \cdot K \left [ \frac{1}{2^k} + \sum_{j=0}^{k-1} \frac{1}{2^k - 2^j} + \sum_{j=k+1}^\infty \frac{1}{2^j - 2^k} \right ] \\
 & \leq & |\lambda| \cdot K \left [ \frac{1}{2^k} + \sum_{j=0}^{k-1} \frac{1}{2^{k-1}} + \sum_{j=k+1}^\infty \frac{1}{2^{j-1}} \right ] \\
 & \leq & |\lambda| \cdot K \left [ \frac{k+2}{2^{k-1}}\right ] \ra 0 \qquad \hbox{as} \qquad k \ra \infty \, .
\end{eqnarray*}
This proves (6.4).   Since the result is true for any unit vector $\xi$, we may conclude that
$\psi$ is a homeomorphism.  Also, if $g$ is any bounded analytic function on the
ball $B$, then there is an $h \in H^\infty(B)$ such that
$$
\widehat{h}(\psi(\lambda)) = g(\lambda) \, .
$$
In fact we have only to take $h = g \circ \widehat{f}$.

In summary, we have constructed a homeomorphism $\psi$ of the open ball $B$ into the fiber ${\cal M}_{\bf 1}$.
Also $\psi$ is analytic, in the sense that $\widehat{g} \circ \psi$ is analytic for
every $g \in H^\infty(B)$.  Therefore the ball $\psi(B)$ has a natural analytic structure.  When we restrict
the algebra $\widehat{H^\infty(B)}$ to this ball, we obtain the algebra of all bounded analytic
function on $\psi(B)$.  It is thus easy to see that the uniformly closed restriction
algebra $\widehat{H^\infty(B)} \bigr |_{\psi(B)}$ will have as its maximal ideal space
the subset $\Xi$ of ${\cal M}$ defined by
$$
\Xi = \{\varphi \in {\cal M}: |\varphi(f)| \leq \sup_{\psi(B)} |\widehat{f}| \ \ \hbox{for all} \ f \in H^\infty(B)\} \, .
$$
This set $\Xi$ is contained in ${\cal M}_{\bf 1}$, as we see by looking at
$f(\lambda) = \bigl ( (1 + \lambda)/2, 0, 0, \dots, 0)$.  Since this restriction
algebra is isomorphic to the algebra of bounded analytic functions on the ball, 
the set $\Xi$ must be homeomorphic to the entire maximal ideal space ${\cal M}$.

The maximum modulus principle now makes it clear that $\psi(B)$ lies in ${\cal M}_{\bf 1} \setminus \Gamma$,
so we see even more clearly that $\overline{\btu} \setminus \btu \ne \Gamma$.  We see from the previous
discussion that the space ${\cal M}$ reproduces itself in any given fiber
{\it ad infinitum}.  Because in $\Xi$ there are fibers attached to the disc $\psi(B)$ in each of which
is a closed set homeomorphic to ${\cal M}$, and so forth.

\section{Gleason Parts}

Although [GLE] is the source of the idea of part, the reference [GAM] has a much more
polished and developed presentation.  We use the latter as the inspiration for our exposition here.
The idea of ``part'' is an interesting equivalence class on the maximal ideal space
${\cal M}$ of $H^\infty$.  

We say that a point $\theta \in {\cal M}$ is in the {\it same part} as $\phi \in {\cal M}$ if there is a
$c > 0$ such that 
$$
\frac{1}{c} < \frac{\widehat{u}(\theta)}{\widehat{u}(\phi)} < c	\eqno (7.1)
$$
for any positive $u$ that is the real part of an element of $H^\infty(B)$.   Gamelin [GAM] refers
to (7.1) as a {\it Harnack inequality}.  The relation expressed by (7.1) is clearly an equivalence
relation, and the equivalence classes are called {\it parts}.

We first note that, if a connected subset $U$ of ${\cal M}$ is to be endowed with an analytic
structure so that the functions in $H^\infty(B)$ become (under the Gelfand transform) analytic on $U$, then of course Harnack's 
inequality will hold on $U$.  Therefore $U$ must lie in a single part of ${\cal M}$.  It is certainly
a matter of some interest to explore to what extent the parts of ${\cal M}$ can be equipped with
analytic structure.

Now assume that $\theta$ and $\phi$ belong to the same part of ${\cal M}$.  Let $b(\theta, \phi)$
be the infimum of all $c$ for which Harnack's inequality is valid for $\theta$ and $\phi$.  Then $b$
has these properties:
\begin{enumerate}
\item[{\bf (a)}]  $\displaystyle \frac{1}{b(\theta, \phi)} \leq \frac{u(\theta)}{u(\phi)} \leq b(\theta, \phi)$ for $u \in \hbox{Re}\, (H^\infty(B))$, $u > 0$;
\item[{\bf (b)}]  $b(\theta, \phi) \geq 1$;
\item[{\bf (c)}]  $b(\theta, \phi) = 1$ if and only if $\theta = \phi$;
\item[{\bf (d)}]  $b(\theta, \phi) = b(\phi, \theta)$;
\item[{\bf (e)}]  $b(\theta, \phi) b(\phi,\psi) \geq b(\theta, \psi)$ if $\psi$ also belongs to the same part as $\theta$.
\end{enumerate}

\noindent In particular, by property {\bf (e)}, we see that $\log b(\theta, \phi)$ is a metric on each part of ${\cal M}$.

These five properties are self-evident, and in fact hold for any function algebra (not just
$H^\infty$).  So we shall not discuss their proofs.

We now want to build up towards A. M. Gleason's original definition of part.   To this end we
have

\begin{theorem} \sl
Let $\theta$ and $\phi$ belong to the same part of ${\cal M}$.  Then there are mutually absolutely
continuous representing measures $\mu$ for $\theta$ and $\nu$ for $\phi$ such that
$$
\frac{1}{b(\theta, \phi)} \leq \frac{d\mu}{d\nu} \leq b(\theta, \phi) \, .
$$
\end{theorem}
We refer the reader to [GAM, p.\ 31] for the concept of representing measure.
\smallskip \\

\noindent {\bf Proof of the Theorem:}  Write, in short, $b$ for $b(\theta, \phi)$.  Since 
$$
b u(\phi) - u(\theta) \geq 0 \, ,
$$ 
we know that, for all positive $u \in \hbox{Re}\, (H^\infty(B))$, there is a positive measure
$\alpha$ on ${\cal M}$ such that
$$
bu(\phi) - u(\theta) = \int u \, d\alpha
$$
for all $u \in \hbox{Re}(H^\infty(B))$.  
Likewise there is a positive measure $\beta$ on ${\cal M}$ such that
$$
b u(\theta) - u(\phi) = \int u \, d\beta
$$
for all $u \in \hbox{Re}(H^\infty(B))$.

Solving these equations for $u(\theta)$ and $u(\phi)$ we see that
\begin{eqnarray*}
\mu & = & \frac{b\beta + \alpha}{b^2 - 1} \\
\nu & = & \frac{b\alpha + \beta}{b^2 - 1}
\end{eqnarray*}
are representing measures for $\theta$ and $\phi$ respectively.  These
measures have the desired properties.
\endpf
\smallskip \\

\begin{corollary} \sl
If $\theta$ and $\phi$ lie in the same part of ${\cal M}$, and if $\eta$
is a representing measure for $\theta$, then there is a representing
measure $\lambda$  for $\phi$ such that $\eta$ is absolutely continuous
with respect to $\lambda$, and furthermore $d\eta/d\lambda \leq b(\theta,\phi)$.
\end{corollary}
{\bf Proof:}  The measure $\lambda = \eta/b + \nu - \mu/b$ does the job.
\endpf
\smallskip \\

\begin{theorem} \sl
The following are equivalent for $\theta, \phi \in {\cal M}$:
\begin{enumerate}
\item[{\bf (i)}]  $\theta$ and $\phi$ are in the same part of ${\cal M}$;
\item[{\bf (ii)}]  $\|\theta - \phi\| < 2$, the norm being that of the 
dual of $H^\infty(B)$;
\item[{\bf (iii)}]  The norm of the restriction of $\theta$ to the kernel of 
$\phi$ is less than 1;
\item[{\bf (iv)}]  Whenever $\{f_j\}$ is a sequence in $H^\infty(B)$ such that $\|f_j\| \leq 1$
and $|f_j(\theta)| \ra 1$, then $|f_j(\phi)| \ra 1$.
\end{enumerate}
\end{theorem}
{\bf Proof:}  Gamelin [GAM, p.\ 144] proves this result for a general uniform
algebra $A$.  So we need not repeat the details here.
\endpf
\smallskip

We note that part {\bf (ii)} of the last theorem is Gleason's original definition of ``part.''

\section{Concluding Remarks}

We are still a long way from solving the corona problem in several complex
variables.  We hope that, in analogy with the papers [SCH] and [GLE], the present
work will help to lay the foundations for an attack on this important problem.

\newpage

\noindent {\Large \sc References}
\vspace*{.2in}

\begin{enumerate}

\newfam\msbfam
\font\tenmsb=msbm10  scaled \magstep1 \textfont\msbfam=\tenmsb
\font\sevenmsb=msbm7 scaled \magstep1 \scriptfont\msbfam=\sevenmsb
\font\fivemsb=msbm5  scaled \magstep1 \scriptscriptfont\msbfam=\fivemsb
\def\Bbb{\fam\msbfam \tenmsb}

\def\RR{{\Bbb R}}
\def\CC{{\Bbb C}}
\def\QQ{{\Bbb Q}}
\def\NN{{\Bbb N}}
\def\ZZ{{\Bbb Z}}
\def\II{{\Bbb I}}
\def\TT{{\Bbb T}}
\def\BB{{\Bbb B}}

\item[{\bf [ALE]}] A. B. Aleksandrov, The existence of inner
functions in the ball, {\it Math. USSR Sbornik} 46(1983),
143--159.

\item[{\bf [CAR]}] L. Carleson, Interpolation by bounded
analytic functions and the corona problem, {\it Ann.\ of
Math.} 76(1962), 542--559.

\item[{\bf [CHG]}] M. Christ and D. Geller, Singular integral
characterizations of Hardy spaces on homogeneous groups, {\it
Duke Math.\ J.} 51(1984), 547--598.

\item[{\bf [GAM]}] T. W. Gamelin, {\it Uniform Algebras},
Prentice-Hall, Englewood Cliffs, NJ, 1969.

\item[{\bf [GEL]}] I. M. Gelfand, Normierte Ringe, {\it Mat.\
Sbornik} N.S. 9(1941), 3--24.

\item[{\bf [GLE]}] A. M. Gleason, Function algebras, in {\it
Seminar on Analytic Functions}, vol.\ II, Princeton, 1957,
pp.\ 213--226.

\item[{\bf [HAS]}] M. Hakim and N. Sibony, Fonctions
holomorphes born\'{e}es sur la boule unit\'{e} de $\CC^n$,
{\it Invent.\ Math.} 67(1982), 213--222.

\item[{\bf [HOF]}] K. Hoffman, Bounded analytic functions and
Gleason parts, {\it Ann. of Math.} 86(1967), 74--111.

\item[{\bf [KRA1]}] S. G. Krantz {\it Function Theory of
Several Complex Variables}, 2nd ed., Chelsea Publishing, the
American Mathematical Society, Providence, RI, 2001.

\item[{\bf [LOO]}] L. Loomis, {\it An Introduction to Abstract
Harmonic Analysis}, Van Nostrand, New York, 1953.

\item[{\bf [LOW]}] E. L\o w, A construction of inner functions
on the unit ball in $\CC^p$, {\it Invent.\ Math.} 67(1982),
223--229.

\item[{\bf [NEV]}] R. Nevanlinna, {\it Eindeutige Analytische
Funktionen} 2 Aufl., Springer,
Berlin-G\"{o}ttingen-Heidelberg, 1953.

\item[{\bf [RUD1]}] W. Rudin, {\it New Constructions of
Functions Holomorphic in the Unit Ball of \boldmath $\CC^n$},
CBMS, American Mathematical Society, Providence, RI, 1986.

\item[{\bf [RUD2]}] W. Rudin, {\it Function Theory in the Unit
Ball of \boldmath $\CC^n$}, Springer-Verlag, New York, 1980.

\item[{\bf [SCH]}] I. J. Schark, Maximal ideals in an algebra
of bounded analytic functions, {\it J. Math.\ Mech.} 10(1961),
735--746.

\end{enumerate}
\vspace*{.25in}

\begin{quote}
Steven G. Krantz  \\
Department of Mathematics \\
Washington University in St. Louis \\
St.\ Louis, Missouri 63130  \\
{\tt sk@math.wustl.edu}
\end{quote}

\end{document}